
\input amstex
\documentstyle{amsppt}
\magnification=\magstep1
\hcorrection{0in}
\vcorrection{0in}
\pagewidth{6.3truein}
\pageheight{9.3truein}
\NoRunningHeads

\def\Z{\Bbb Z}
\def\Q{\Bbb Q}
\def\P{\Bbb P}
\def\O{\Cal O}
\def\I{\Cal I}
\def\Bs{\text{Bs}}
\def\Supp{\text{Supp }}
\def\Div{\text{Div }}
\def\ccdot{\cdot\cdot}
\topmatter
\title On freeness theorem of the adjoint bundle on a normal surface
\endtitle
\author Takeshi Kawachi \endauthor
\affil Department of Mathematics, Tokyo Institute of Technology \endaffil
\address 2-12-1 Oh-okayama, Meguro-ku, Tokyo, JAPAN, 152 \endaddress
\email kawachi\@math.titech.ac.jp \endemail
\abstract
The adjoint linear system on a surface have been studied by many authors.
Among these, Reider's criterion on a smooth surface is very famous.
Here we prove a similar criterion on normal surfaces.
\endabstract
\endtopmatter
\document
\head Introduction \endhead
Let $k$ be an algebraically closed field of characteristic $0$.
Throughout this paper, we consider everything over $k$.\par
It is an interesting problem when the adjoint linear system
$|K_X+mL|$ is free.
Reider proved the following theorem in [R].\par
\proclaim{Theorem}
Let $S$ be a smooth complex algebraic surface and $L$ be a nef divisor
on $S$.
If $L^2\geq 5$ and $p$ is a base point of $|K_X+L|$,
then there exists an effective divisor $E$ passing through $p$
such that
$$\matrix
  \format \r & \quad \l \\
  \text{either} & LE=0,\ E^2=-1 \\
  \text{or} & LE=1,\ E^2=0. \\
  \endmatrix
$$
\endproclaim\par
In [S2], Sakai extended this theorem to normal surfaces as follows.\par
\proclaim{Theorem}
Let $Y$ be a normal Moishezon surface.
Let $D$ be a nef divisor on $Y$ such that $K_Y+D$ is a Cartier divisor.
If $D^2>4+\eta(Y)$, then $|K_Y+D|$ has no base points
unless there exists a nonzero effective divisor $E$ such that
$$\gather
  0\leq DE<2+\frac12\eta(Y) \cr
  DE-1-\frac14\eta(Y)\leq E^2\leq\frac{(DE)^2}{D^2} \cr
  E^2<0 \text{ if } DE=0.
\endgather$$
\endproclaim\par
In this theorem, the number $\eta(Y)$ is a non-negative rational number
determined by the singularities of $Y$
which may be arbitrarily large if $Y$ has bad singularities.
On the other hand, Matsushita
proved the following theorem in [Ma].\par
\proclaim{Theorem}
Let $Y$ be a normal projective surface
which has only log terminal singularities.
Fix one point $x_0$ on $Y$.
Let $D$ be a nef big $\Q$-divisor which satisfies following conditions.
\roster
\item $D^2>4$.
\item $D.C>2$ for all curves $C\subset Y$ such that $x_0\in C$.
\item $K_Y+\lceil D\rceil$ is a Cartier divisor.
\endroster
Then $x_0\not\in\Bs |K_Y+\lceil D\rceil|$.
\endproclaim\par
We will generalize and improve the above results as follows.
First we introduce the terminology.
\definition{Definition}
Let $Y$ be a normal surface.
$Y$ is said to have log terminal singularity
if the following conditions are satisfied:\newline
\quad (\hbox to 2ex{\hfil i\hfil})
$K_Y\in\Div(Y)\otimes\Q$.\newline
\quad (\hbox to 2ex{\hfil ii\hfil})
There exists a resolution of singularities
$f\: X\to Y$ such that the union of the exceptional locus
is a divisor with only normal crossings, and that
$$
  K_X=f^*K_Y+\sum a_j\Delta_j\text{ for $a_j\in\Q$ }
$$
with the condition that $a_j>-1$ whenever $\Delta_j$ is exceptional for $f$.
\enddefinition\par
Let $(Y,x)$ be a germ of normal surface singularity.
Let $f\:X\to Y$ be the minimal resolution.
There exists an effective $\Q$-divisor $\Delta$ supported on
$f^{-1}(x)$ such that $f^*K_Y=K_X+\Delta$
($\Delta$ has the opposite sign of Kawamata log discrepancy).
Let $Z$ be the fundamental cycle of $x$.\par
If $x$ is a log terminal singular point then we obtain the following.\par
\proclaim{Proposition 1}
Let $Y,x,\Delta,Z$ as above.
If $x$ is a log terminal singular point of $Y$ then $-(\Delta-Z)^2\leq 2$.
Furthermore $-(\Delta-Z)^2=2$ if and only if $x$ is a rational double point.
\endproclaim\par
We define
$$
  \delta_x=
    \cases
      4, & \text{$x$ is a smooth point,} \\
      -(\Delta-Z)^2, & \text{$x$ is a log terminal singular point,} \\
      0, & \text{$x$ is a singular point which is not log terminal.}
    \endcases
$$
The main result of this paper is the following.\par
\proclaim{Theorem}
Let $Y$ be a normal surface over $k$.
Let $D$ be a nef divisor on $Y$ such that $K_Y+D$ is a Cartier divisor.
Let $x$ be a point on $Y$ and $D^2>\delta_x$.
If $x$ is a base point of $|K_Y+D|$ then
there exists a non-zero effective divisor $E$ on $Y$
passing through $x$ such that
$$\gather
  0\leq DE<\frac12\delta_x \cr
  DE-\frac14\delta_x\leq E^2\leq\frac{(DE)^2}{D^2} \cr
  E^2<0 \text{ if } DE=0. \cr
\endgather$$
\endproclaim\par
In particular, if $x$ is a singular point which is not log terminal
then $|K_Y+D|$ is always free at $x$.
One may say that, the worse singular the point $x$ is,
the less likely it is a base point of the adjoint system.
By this theorem, we obtain the following results immediately.\par
\proclaim{Corollary 1}
Let $Y$ be a normal surface over $k$.
Let $A$ be an ample divisor on $Y$.
If $t\geq 3$ and $K_Y+tA$ is a Cartier divisor then $|K_Y+tA|$ is free.
Furthermore, if $A^2>1$ and $K_Y+2A$ is a Cartier divisor
then $|K_Y+2A|$ is free.
\endproclaim\par
This is one of the normal surface case of Fujita conjecture.\par
\proclaim{Corollary 2}
Let $Y$ be a normal surface over $k$.
Let $D$ be a nef divisor on $Y$.
Fix one point $x$ in $Y$.
If $D$ satisfies
\roster
\item $D^2>\delta_x$
\item $DC\geq \frac12\delta_x$ for all irreducible curves $C$ on $Y$
\item $K_Y+D$ is a Cartier divisor
\endroster
then $|K_Y+D|$ is free at $x$.
\endproclaim\par
\proclaim{Corollary 3}
Let $Y$ be a normal surface having only rational double points as
singularities.
If $D$ is an ample Cartier divisor and $D^2>2$ then
every base points of $|K_Y+D|$ is smooth.
\endproclaim\par
We prove Proposition 1 in \S1.
We prove the main theorem in \S2.
The proof of the theorem is divided to three cases.
First we treat the case that $x$ is a non log terminal singular point.
Second we treat the case that $x$ is a smooth point.
Finally we treat the case that $x$ is a log terminal singular point.
\par
\subhead Acknowledgement\endsubhead\par
The author would like to thank to Professor T.~Fujita for
his useful advice.
By his advice, the author can improve the main theorem
and its proof.\par
\head \S1 Preliminaries \endhead
Let $(Y,x)$ be a germ of normal surface singularity.
Let $f\:X\to Y$ be the minimal resolution of $x$.
There exists an effective $\Q$-divisor $\Delta$ supported on
$f^{-1}(x)$ such that $f^*K_Y=K_X+\Delta$.
Let $Z$ be the fundamental cycle of $x$.\par
\proclaim{Proposition 1}
Let $Y,x,\Delta,Z$ as above.
If $x$ is a log terminal singular point of $Y$ then $-(\Delta-Z)^2\leq 2$.
Furthermore $-(\Delta-Z)^2=2$ if and only if
$x$ is a rational double point.
\endproclaim\par
\demo{Proof}
If $x$ is a rational double point then $\Delta=0$ and $-Z^2=2$.
Hence the equation $-(\Delta-Z)^2=2$ holds.\par
Let $x$ be a log terminal singular point but not a rational double point.
Let $\Delta=\sum_{i=1}^n a_i C_i$ where $0< a_i<1$ for all $i$
where $C_n$ and $C_{n-1}$ etc.\ in the case of type $D_n$ or $E_n$ mean
$C_1''$ and $C_1'$ etc.\ respectively.
In such cases, $a_n$ etc.\ are written as $a_1''$ etc.
The rational numbers $a_i$ are defined by the equations
$K_X C_i+\Delta C_i=0$ for all $i$.
Since $(K_X+C_i)C_i=-2$, we have $\Delta C_i=-K_X C_i=C_i^2+2$.\par
we divide the cases by the type of the singularity (cf\. Appendix).\par
Case 1: The dual graph is of type $A_n$. \par
In this case, the fundamental cycle $Z$ is $Z=\sum C_i$.
If $n=1$ then $a_1C_1^2=C_1^2+2$.
Hence $(\Delta-Z)C_1=2$.
Therefore we have $-(\Delta-Z)^2=2(1-a_1)=-4/C_1^2$.
Now we assume $n\geq 2$.
Hence we have the equations
$$\align
  & (\Delta-Z)C_1 = C_1^2+2-(C_1^2+1) = 1 \\
  & (\Delta-Z)C_i = C_i^2+2-(C_i^2+2) = 0\ (2\leq i\leq n-1) \\
  & (\Delta-Z)C_n = C_n^2+2-(C_n^2+1) = 1.
\endalign$$
Therefore we have
$$
  -(\Delta-Z)^2 = -(a_1-1)-(a_n-1) = 2-a_1-a_n < 2.
$$
\par
Case 2: The dual graph is of type $D_n$.\par
Let $Z=\sum_{i=1}^n b_i C_i$ be the fundamental cycle of $x$.
If all $C_i^2=-2$ then $x$ is a rational double point.
Hence there exists some $j$ such that $C_j^2\leq -3$.\par
($D_n$-i) First, let $C_{n-2}^2\leq -3$.
Then the fundamental cycle is $\sum_{i=1}^{n-2} C_i +C'_1+C''_1$.
Hence we have
$$\align
  & (\Delta-Z)C_1 = C_1^2+2-(C_1^2+1) = 1 \\
  & (\Delta-Z)C_i = C_i^2+2-(C_i^2+2) = 0\ (2\leq i\leq n-3) \\
  & (\Delta-Z)C_{n-2} = C_{n-2}^2+2-(C_{n-2}^2+3) = -1 \\
  & (\Delta-Z)C'_1 = {C'_i}^2+2-({C'_1}^2+1) = 1 \\
  & (\Delta-Z)C''_1 = {C''_i}^2+2-({C''_1}^2+1) = 1.
\endalign$$
Therefore we have
$$\align
  -(\Delta-Z)^2
  &= -(a_1-1)+(a_{n-2}-1)-(a'_1-1)-(a''_1-1) \\
  &= 2-a_1+a_{n-2}-a'_1-a''_1.
\endalign$$
Since $\Delta C'_1=\Delta C''_1=0$,
the numbers $a_{n-2},\ a'_1,\ a''_1$ satisfy the equations
$$
  a_{n-2}-2a'_1=0,\ a_{n-2}-2a''_1=0.
$$
Hence we have $a_{n-2}-a'_1-a''_1=0$.
Therefore $-(\Delta-Z)^2=2-a_1< 2$. \par
($D_n$-ii) Next, let $C_{j+1}^2=\cdots =C_{n-2}^2=-2$ and $C_j^2\leq -3$.
Then the fundamental cycle $Z$ is
$\sum_{i=1}^j C_i+\sum_{i=j+1}^{n-2} 2C_i+C'_1+C''_1$.
If $2\leq j\leq n-3$ then we have the equations
$$\align
  & (\Delta-Z)C_1 = C_1^2+2-(C_1^2+1) = 1 \\
  & (\Delta-Z)C_i = C_i^2+2-(C_i^2+2) = 0\ (2\leq i\leq j-1) \\
  & (\Delta-Z)C_j = C_j^2+2-(C_j^2+3) = -1 \\
  & (\Delta-Z)C_{j+1} = C_{j+1}^2+2-(2C_{j+1}^2+3) = 1 \\
  & (\Delta-Z)C_i = C_i^2+2-(2C_i^2+4) = 0\ (j+2\leq i\leq n-2) \\
  & (\Delta-Z)C'_1 = {C'_1}^2+2-({C'_1}^2+2) = 0 \\
  & (\Delta-Z)C''_1 = {C''_1}^2+2-({C''_1}^2+2) = 0.
\endalign$$
Therefore we have
$$
  -(\Delta-Z)^2=-(a_1-1)+(a_j-1)-(a_{j+1}-2)=2-a_1+a_j-a_{j+1}.
$$
In this case, the numbers $a_j,\dots ,a_{n-2},a'_1,a''_1$
satisfy the equations
$$\align
   & a_j-2a_{j+1}+a_{j+2}=0,\dots ,a_{n-4}-2a_{n-3}+a_{n-2}=0, \\
   & a_{n-2}=2a'_1,\ a_{n-2}=2a'_1, \\
   & a_{n-3}-2a_{n-2}+a'_1+a''_1=a_{n-3}-a_{n-2}=0.
\endalign$$
Hence we have $ a_j=a_{j+1}=\cdots =a_{n-2}$.
Therefore $-(\Delta-Z)^2=2-a_1< 2$.\par
($D_n$-iii) Let $j=1$.
Then we have the equations
$$\align
  & (\Delta-Z)C_1 = C_1^2+2-(C_1^2+2) = 0 \\
  & (\Delta-Z)C_{j+1} = C_2^2+2-(2C_2^2+3) = 1 \\
  & (\Delta-Z)C_i = C_i^2+2-2(C_i^2+2) = 0\ (3\leq i\leq n-2) \\
  & (\Delta-Z)C'_1 = {C'_1}^2+2-({C'_1}^2+2) = 0 \\
  & (\Delta-Z)C''_1 = {C''_1}^2+2-({C''_1}^2+2) = 0.
\endalign$$
Hence we have $-(\Delta-Z)^2=2-a_2< 2$.\par
Case 3: The dual graph is of type $E_n$.\par
If the dual graph is in the Appendix then we write $(\mu;a,b,c;d,e)$.
We set indexes $i$ as in the Appendix.\par
($E_n$-i)First, we assume $\mu\geq 3$.
Then the fundamental cycle is $\sum C_i+\sum C'_i+C''_1$.
In this case, we have the equations
$$\align
  & (\Delta-Z)C_1 = C_1^2+2-(C_1^2+1) = 1 \\
  & (\Delta-Z)C_i = C_i^2+2-(C_i^2+2) = 0\ (2\leq i\leq n_1) \\
  & (\Delta-Z)C'_1 = {C'_1}^2+2-({C'_1}^2+1) = 1 \\
  & (\Delta-Z)C'_i = {C'_i}^2+2-({C'_i}^2+2) = 0\ (2\leq i\leq n_2) \\
  & (\Delta-Z)C''_1 = {C''_1}^2+2-({C''_1}^2+1) = 1 \\
  & (\Delta-Z)C_0 = C_0^2+2-(C_0^2+3) = -1. \\
\endalign$$
Therefore we have
$$\align
  -(\Delta-Z)^2
  &= -(a_1-1)-(a_{n-2}-1)-(a_{n-1}-1)+(a_n-1) \\
  &= 2-a_1-a'_1-a''_1+a_0.
\endalign$$
Since the numbers $a_0,\ a''_1$ satisfy the equation $a_0-2a''_1=0$,
we have
$$
  -(\Delta-Z)^2=2-a_1-a'_1+a''_1.
$$
Now we calculate $\delta=2-a_1-a'_1+a''_1$ for all cases of type $E_n$.
The results are listed in Table 1 in the Appendix.
Hence if $\mu\geq 3$ then $\delta<2$. \par
($E_n$-ii) Next, we assume $\mu=2$.
Then the fundamental cycle is as in the Appendix.
We can calculate $\delta=-(\Delta-Z)^2$ directly,
we have the following result.
$$\vbox{\offinterlineskip
  \halign{&\vrule#&\strut\ \hfil#\hfil\vbox{\hrule width0pt height12pt} \cr
    \noalign{\hrule}
    & type && 2 && 3 && 5 && 6 && 7 && 9 && 10 &&
      11 && 12 && 13 && 14 && 15 & \cr
    height2pt & \omit && \omit && \omit && \omit && \omit && \omit &&
    \omit && \omit && \omit && \omit && \omit && \omit && \omit & \cr
    \noalign{\hrule}
    & $\delta$ && $\frac{14}{9}$ && $\frac65$ && $\frac85$ &&
      $\frac{10}{7}$ && $\frac{12}{11}$ && $\frac{18}{11}$ &&
      $\frac{23}{13}$ && $\frac{33}{23}$ && $\frac{11}{7}$ &&
      $\frac{21}{17}$ && $\frac{26}{19}$ && $\frac{30}{29}$ & \cr
    height2pt & \omit && \omit && \omit && \omit && \omit && \omit &&
    \omit && \omit && \omit && \omit && \omit && \omit && \omit & \cr
    \noalign{\hrule}
  }}
$$
Hence we have $-(\Delta-Z)^2< 2$. \qed
\enddemo\par
\head \S2 Main result \endhead
Let $D$ be any $\Q$-divisor.
We write the integral part of $D$ as $[D]$,
and the fractional part of $D$ as $\{D\}$.
Let $Y$ be a normal surface and $D$ be a $\Q$-divisor on $Y$.
Let $f\:X\to Y$ be the minimal resolution of singularities on $Y$.
Let $\Delta=f^*K_Y-K_X$.
If $M=K_Y+D$ is a Cartier divisor then $f^*M$ is an integral divisor.
Hence $f^*M-K_X=f^*D+\Delta$ is an integral divisor.
Therefore $D$ turns out to be a $\Z$-coefficient $\Q$-Cartier divisor.\par
Let $x$ be a singular point on $Y$ and $Z$ be the fundamental cycle of $x$.
Let $\Delta_x$ be the components of $\Delta$ supported on $f^{-1}(x)$.
We define
$$
  \delta_x=
    \cases
      4, & \text{$x$ is a smooth point,} \\
      -(\Delta_x-Z)^2, & \text{$x$ is a log terminal singular point,} \\
      0, & \text{$x$ is a singular point which is not log terminal.}
    \endcases
$$
\proclaim{Theorem}
Let $Y$ be a normal surface over $k$.
Let $D$ be a nef divisor on $Y$ such that $K_Y+D$ is a Cartier divisor.
Let $x$ be a point on $Y$ and $D^2>\delta_x$.
If $x$ is a base point of $|K_Y+D|$ then
there exists a non-zero effective divisor $E$ on $Y$
passing through $x$ such that
$$\gather
  0\leq DE<\frac12\delta_x \cr
  DE-\frac14\delta_x\leq E^2\leq\frac{(DE)^2}{D^2} \cr
  E^2<0 \text{ if } DE=0. \cr
\endgather$$
\endproclaim\par
{\it Proof.}
Let $f,\ \Delta,\ Z$ as above.
We divide the cases by the type of $x$.\par
\subhead 1. The case where $x$ is not a log terminal point:
\endsubhead\par
Let $\Delta_x=\sum a_i \Delta_i$ be the prime decomposition.\par
\proclaim{Lemma 1}
$a_i\geq 1$ for some $i$.
\endproclaim\par
\demo{Proof}
Suppose that $0<a_i<1$ for all $i$.
Since $K_X\Delta_i=-\Delta\Delta_i$,
we have
$$
  (K_X+\Delta_i)\Delta_i
  = -\sum_{j\ne i} a_j\Delta_j\Delta_i+(1-a_i)\Delta_i^2<0.
$$
Hence each $\Delta_i\cong\P^1$.\par
We claim $\Delta_i\Delta_j\leq 1$ for all $i,j$.
Indeed, if $\Delta_i\Delta_j\geq 2$ for some $i,j$ with $i\ne j$ then
$$\align
  -2 &= (K_X+\Delta_i)\Delta_i\leq (1-a_i)\Delta_i^2-a_j\Delta_i\Delta_j\\
     &\leq -2(1-a_i)-2a_j.
\endalign$$
Thus we have $a_j\leq a_i$.
Similarly, if we exchange the roles of $i$ and $j$,
we have $a_i\leq a_j$.
Hence $a_i=a_j$ and all inequality must be equality.
Thus $\Delta_i^2=\Delta_j^2=-2$ and $\Delta_i\Delta_j=2$.
But this shows $(\Delta_i+\Delta_j)^2=0$, which is absurd.\par
Suppose that three or more components meet at one point.
If $\Delta_i$, $\Delta_j$ and $\Delta_k$ meet at a point
and $a_i\leq a_j\leq a_k$ then
$$\align
  -2 &= (K_X+\Delta_i)\Delta_i\leq (1-a_i)\Delta_i^2-a_j-a_k \\
     &\leq -2(1-a_i)-a_j-a_k\leq -2.
\endalign$$
Thus we have $\Delta_i^2=-2$ and $a_i=a_j=a_k$.
If we exchange the roles of $i$, $j$ and $k$,
we have $\Delta_j^2=\Delta_k^2=-2$.
But this shows $(\Delta_i+\Delta_j+\Delta_k)^2=0$, which is absurd.\par
Hence $\Delta_x$ is a normal crossing divisor
and all coefficients are less than $1$.
Thus $x$ is a log terminal singularity, which contradicts our assumption.
\qed
\enddemo\par
By the above lemma we have $[\Delta_x]\ne 0$.
Let $M=K_Y+D$.
Then $f^*M-K_X-\Delta=f^*D$ is nef and big.
Hence we have
$$
   H^1(X,\O(f^*M-[\Delta]))=H^1(X,\O(K_X+f^*D+\{\Delta\}))=0
$$
from the following vanishing theorem due to Miyaoka.\par
\proclaim{Theorem [Mi]}
Let $D$ be a divisor on a surface $X$ with $\kappa (X,D)=2$
and $D=P+\sum n_iN_i$ its Zariski decomposition.
If each $n_i$ is smaller than $1$, then $H^1(X,-D)$ vanishes.
\endproclaim\par
Consider the following exact sequence
$$\align
   H^0(X,\O(f^*M-[\Delta]+[\Delta_x])) &\to
   H^0([\Delta_x],\O(f^*M-[\Delta]+[\Delta_x])|_{[\Delta_x]}) \\
   & \to H^1(X,\O(f^*M-[\Delta]))=0.
\endalign$$
Since $[\Delta]-[\Delta_x]=\sum_{p\ne x} [\Delta_p]$
where $p$ is a singular point on Y,
$[\Delta]-[\Delta_x]$ is disjoint from $[\Delta_x]$.
Hence we have
$$
   H^0([\Delta_x],\O(f^*M-[\Delta]+[\Delta_x])|_{[\Delta_x]})
   \cong H^0([\Delta_x],\O(f^*M)|_{[\Delta_x]}).
$$
Since $f([\Delta_x])=\{x\}$ is a point on $Y$,
we have
$$
   H^0([\Delta_x],\O(f^*M)|_{[\Delta_x]})
   \cong H^0([\Delta_x],\O)\supset k.
$$
Therefore there must be a $\varphi\in H^0(X,\O(f^*M-\sum_{p\ne x}[\Delta_p]))$
such that $\varphi\equiv 1$ on $[\Delta_x]$.
Hence the injection
$$
   H^0(X,\O(f^*M-\sum_{p\ne x} [\Delta_p]))\to
   H^0(X,\O(f^*M))\cong H^0(Y,\O(M))
$$
gives $\psi\in H^0(Y,M)$ such that $\psi(x)\ne 0$.
Therefore $x$ is not a base point of $|M|=|K_Y+D|$.\par
We complete the proof of non log terminal case.\par
\vskip\baselineskip
\subhead 2. The case where $x$ is a smooth point:
\endsubhead\par
In this case $\delta_x=4$.
Let $\pi\:\widetilde{X}\to X$ be the blowing up of $X$ at $x$.
Let $L$ be the exceptional curve of $\pi$.
Since $(\pi^*f^*D-2L)^2=D^2-4>0$ by the assumption and
$\pi^*f^*D(\pi^*f^*D-2L)=D^2>0$ for nef divisor $\pi^*f^*D$,
we have $\pi^*f^*D-2L$ is big.
Hence $\pi^*f^*D-2L$ has the Zariski decomposition
$\pi^*f^*D-2L=P+N$ where $P$ is a nef $\Q$-divisor and
$N$ is negative definite.
Since $(\pi^*f^*D-2L)^2=D^2-4>0$, we have $P^2>0$.
Hence $P$ is a nef and big $\Q$-divisor.\par
Furthermore $N$ does not contain $L$ as a component.
Indeed, if $PL=0$ then $\pi^*f^*D=P+(N+2L)$ is another Zariski decomposition
of a nef and big divisor $\pi^*f^*D$, which is absurd.
Hence $PL>0$. \par
We have $P\Gamma=0$ for any irreducible $f$-exceptional divisor $\Gamma$,
where we identify $\Gamma$ and $\pi^*\Gamma$.
Indeed, if $P\Gamma>0$ then $N\Gamma<0$ because $(\pi^*f^*D-2L)\Gamma=0$.
Hence $\Gamma$ is a component of $\Supp N$, which is contradiction.\par
Thus we have $P\Delta=0$,
hence $P+(N+\Delta)$ is the Zariski decomposition of $\pi^*f^*D+\Delta-2L$.
Let $\widetilde{D}=\pi^*f^*D+\Delta-2L=\pi^*f^*(K_Y+D)-K_{\widetilde{X}}-L$.
\par
Consider the following diagram
$$\CD
  H^0(\widetilde{X},\pi^*f^*\O(K_Y+D)) @>>>
    H^0(L,\pi^*f^*\O(K_Y+D)\otimes\O/\O(-L)) \\
  @A{\wr}AA \bigcup \\
  H^0(Y,\O(K_Y+D)) @>>> H^0(Y,\O(K_Y+D)\otimes\O/\I_x) \\
\endCD$$
where $\I_x$ is the maximal ideal sheaf of $x$.
Since $x$ is a base point of $K_Y+D$,
the bottom morphism $H^0(K_Y+D)\to H^0(\O(K_Y+D)\otimes\O/\I_x)$
is not surjective.
Hence we have
$$
  H^1(\widetilde{X},\pi^*f^*\O(K_Y+D)\otimes\O(-L))
  = H^1(\widetilde{X},\O(K_{\widetilde{X}}+\widetilde{D}))
  \ne 0.
$$
On the other hand, since $P$ is nef and big, we have
$$
  H^1(\widetilde{X},\O(K_{\widetilde{X}}+P+\{N+\Delta\}))
  =0
$$
by the vanishing theorem due to Miyaoka.\par
Let $D_0=\widetilde{D}-[N+\Delta]$.
Then we have $H^1(\widetilde{X},\O(-D_0))=0$ by duality.
Let $D_j=D_{j-1}+E_j$ for some irreducible component $E_j$ of
$\widetilde{D}-D_{j-1}$.
Then $D_n=\widetilde{D}$ where $n$ is the number of components of
$[N+\Delta]$ counted with multiplicity.
If $D_{j-1}E_j>0$ for all $j$, then the following sequence
$$
  0=H^0(E_j,\O_{E_j}(-D_{j-1}))\to H^1(\widetilde{X},\O(-D_j))
  \to H^1(\widetilde{X},\O(-D_{j-1}))
$$
gives $H^1(\widetilde{X},\O(-D_n))=0$.
But this contradicts
$H^1(\widetilde{X},\O(K_{\widetilde{X}}+\widetilde{D}))\ne 0.$
Hence there exists some $j$ such that
$H^1(\widetilde{X},\O(-D_j))=0$ and $D_jC\leq 0$
for any irreducible component $C$ of $\widetilde{D}-D_j$.
Let $E=\widetilde{D}-D_j$.
Then $E$ is a non-zero effective divisor and
$E\subset [N+\Delta]$ by our construction.\par
Since $H^1(\widetilde{X},\O(-D_j))=0$, we have
$$
  H^1(\widetilde{X},\O(K_{\widetilde{X}}+\pi^*f^*D+\Delta-2L-E))=0.
$$
Suppose $E$ is disjoint from $L$.
Then $E=\pi^*\pi_* E$ and
$$
  H^1(X,\O(K_X+f^*D+\Delta-\pi_*E)\otimes\I_x)=0.
$$
Consider the following diagram
$$\CD
  H^0(X,\O(K_X+f^*D+\Delta-\pi_*E)) @>>>
    H^0(X,\O(K_X+f^*D+\Delta-\pi_*E)\otimes\O/\I_x) \\
  @VV{+\pi_*E}V @V{\wr}VV \\
  H^0(X,\O(f^*(K_Y+D))) @>>> H^0(X,\O(f^*(K_Y+D))\otimes\O/\I_x). \\
\endCD$$
This diagram shows that $x$ is not a base point of $|K_Y+D|$,
which contradicts our assumption.
Hence $E$ meets $L$.\par
Since $N+\Delta$ does not contain $L$ as a component and
$E$ is taken from $[N+\Delta]$,
this $E$ has a component which is not $\pi$-exceptional and meets $L$.
Therefore $f_*\pi_*E\ne 0$.\par
Let $E-\Delta=E_1-\Delta_1$ where $E_1$ and $\Delta_1$ are
effective $\Q$-divisors without common components.
Since $D_jC\leq 0$ for any component $C$ in $E$, we have
$$\align
  0 &\geq (\pi^*f^*D+\Delta-2L-E)E_1 \\
    &= (\pi^*f^*D-2L-E_1+\Delta_1)E_1 \\
    &\geq (\pi^*f^*D-2L-E_1)E_1.
\endalign$$
Since $E\cap L\ne\emptyset$, we have also $E_1\cap L\ne\emptyset$.
Furthermore since $E$ does not contain $L$ as a component,
$E_1$ does not contain $L$ as a component.
\par
Let $E_2=\pi_*E_1$ and $E_1=\pi^*E_2-aL$.
Then we have
$$\align
  0 &\geq (\pi^*f^*D-2L-\pi^*E_2+aL)(\pi^*E_2-aL) \\
    &= (f^*D-E_2)E_2 +a(a-2) \\
    &\geq (f^*D-E_2)E_2-1.
\endalign$$
Since $E_1$ is not $f$-exceptional, we have $E_3=f_*E_2\ne\emptyset$.
Let $E_2=f^*E_3-\Gamma$.
Then we have
$$\align
  0 &\geq (f^*D-f^*E_3+\Gamma)(f^*E_3-\Gamma)-1 \\
    &= (D-E_3)E_3 -\Gamma^2-1 \\
    &\geq (D-E_3)E_3 -1.
\endalign$$\par
On the other hand,
since $(\pi^*f^*D-2L-E_1)E_1\leq 0$, we have
$$
  DE_3=\pi^*f^*D E_1 < 2LE_1=2a.
$$
Suppose $a\geq2$.
Since $(N+\Delta)-E_1$ is effective and does not contain $L$ as a component,
we have $NL=(N+\Delta)L\geq E_1L=a\geq 2$.
Since $PL>0$, we get
$2=(\pi^*f^*D-2L)L=(P+N)L>2$, which is contradiction.
Thus we infer $a\leq 1$, since $a=E_1L=EL$ and
$E$ and $L$ are $\Z$-divisors on $\widetilde{X}$.
In particular $E_3$ has the component passing through $x$
which is smooth at $x$ because $E$ does not contain $L$ and $EL=1$.
\par
If we replace $E_3$ by E, we have the following.
$$\gather
  0 \leq DE < 2 \\
  DE-1\leq E^2\leq \frac{(DE)^2}{D^2} \\
  E^2 < 0 \text{ if } DE=0.
\endgather$$
Thus we complete the proof when $x$ is a smooth point.\par
\vskip\baselineskip
\subhead 3. The case where $x$ is a log terminal point:
\endsubhead\par
Let $\Delta=\Delta_x+\Delta'$ where $\Delta_x$ consists of the components
supported on $f^{-1}(x)$ and $\Delta'$ is the others.
Let $Z$ be the fundamental cycle of $x$.
In this case $\delta_x=-(\Delta_x-Z)^2\leq 2$ by Proposition 1.
We have $(f^*D+\Delta_x-Z)^2=D^2+(\Delta_x-Z)^2> 0$
by the assumption.
Hence $f^*D+\Delta_x-Z$ is big, this has the Zariski decomposition
$f^*D+\Delta_x-Z=P+N$ where $P$ is a nef and big $\Q$-divisor
and $N$ is negative definite.
Let $\Gamma$ be an irreducible $f$-exceptional divisor with
$f(\Gamma)\ne\{x\}$.
Since $(f^*D+\Delta_x-Z)\Gamma=0$, we have $P\Gamma=0$ as before.
Therefore $f^*D+\Delta-Z=P+(N+\Delta')$ is the Zariski decomposition
of $f^*D+\Delta-Z$.\par
Consider the following diagram
$$\CD
  H^0(X,f^*\O(K_Y+D)) @>>> H^0(X,f^*\O(K_Y+D)\otimes\O/\O(-Z)) \\
  @A{\wr}AA \bigcup \\
  H^0(Y,\O(K_Y+D)) @>>> H^0(Y,\O(K_Y+D)\otimes\O/f_*\O(-Z)). \\
\endCD$$
Since $x$ is a base point of $|K_Y+D|$, we have
$$
  H^1(X,f^*\O(K_Y+D)\otimes\O(-Z))
  =H^1(X,\O(K_X+f^*D+\Delta-Z))\ne 0.
$$
On the other hand, since $P$ is nef and big, we have
$$
  H^1(X,\O(K_X+P+\{N+\Delta'\}))=0
$$
by the vanishing theorem due to Miyaoka.
By the same argument as smooth case,
we have a nonzero effective divisor $E\subset [N+\Delta']$ such that
$$\gather
  H^1(X,\O(K_X+f^*D+\Delta-Z-E))=0 \\
  (f^*D+\Delta-Z-E)C\leq 0\quad
    \text{ for any irreducible component $C$ in $E$.}
\endgather$$
If $E$ is disjoint from $Z$ then the following diagram
$$\CD
  H^0(X,\O(K_X+f^*D+\Delta-E)) @>>>
    H^0(X,\O(K_X+f^*D+\Delta-E)\otimes\O/\O(-Z)) \\
  @VV{+E}V @V{\wr}VV \\
  H^0(X,\O(K_X+f^*D+\Delta)) @>>>
    H^0(X,\O(K_X+f^*D+\Delta)\otimes\O/\O(-Z)) \\
\endCD$$
shows that the bottom morphism is surjective.
Since $x$ is a base point of $|K_Y+D|$, this is contradiction.
Hence $E$ meets $Z$.\par
Let $E-\Delta'=E_1-\Delta_1'$ where $E_1$ and $\Delta_1'$ are
effective $\Q$-divisors without common components.
Thus we have
$$\align
  0 &\geq (f^*D+\Delta-Z-E)E_1 \\
    &= (f^*D+\Delta_x-Z-E_1+\Delta_1')E_1 \\
    &\geq (f^*D+\Delta_x-Z-E_1)E_1.
\endalign$$
Since $E$ meets $Z$ we have $E_1\cap Z\ne\emptyset$.
Note that the coefficient of a component of $E_1$ which meets $Z$ is integral.
\par
We claim that $f_*E_1$ has a component passing $x$.
Indeed, if not, $E_1=E_1'+\Gamma$
where $E_1'$ consists of the components supported on $f^{-1}(x)$
and $\Gamma$ is disjoint from $Z$.
Since $E_1$ is integral near $Z$, $E_1'$ is a divisor.
Thus $\Delta''=\Gamma-\Delta_1'+\Delta'=E-E_1'$ is an effective $\Z$-divisor
and
disjoint from $f^{-1}(x)$.
Since $H^1(X,\O(K_X+f^*D+\Delta-Z-E))=0$,
we have
$$
  H^1(X,\O(K_X+f^*D+\Delta-Z-E_1'-\Delta''))=0.
$$
Consider the following diagram
$$\CD
  H^0(X,\O(K_X+f^*D+\Delta-\Delta'')) @>>>
    H^0(Z+E_1',\O(K_X+f^*D+\Delta-\Delta'')|_{Z+E_1'}) \\
  @VV{+\Delta''}V @V{\wr}VV \\
  H^0(X,f^*\O(K_Y+D)) @>>> H^0(Z+E_1',f^*\O(K_Y+D)|_{Z+E_1'}).
\endCD$$
Hence the bottom morphism is surjective.
Since $x$ is a base point of $|K_Y+D|$, this is contradiction.
Thus $f_*E_1$ is passing through $x$.\par
Let $E_2=f_*E_1$ and $E_1=f^*E_2-\Gamma$.
Then we have the inequality
$$\align
  0 &\geq (f^*D+\Delta_x-Z-f^*E_2+\Gamma)(f^*E_2-\Gamma) \\
    &= (D-E_2)E_2-\left(\Gamma+\frac{\Delta_x-Z}{2}\right)^2
       +\frac{(\Delta_x-Z)^2}{4} \\
    &\geq (D-E_2)E_2-\frac14\delta_x.
\endalign$$\par
Next, we prove the inequality $DE_2<\frac12\delta_x$.\par
First we prove the following two lemmas.\par
\proclaim{Lemma 2}
Let $Z'$ be an irreducible $f$-exceptional divisor.
If $(Z-\Delta_x)Z'\geq 0$ then $PZ'=0$.
\endproclaim\par
\demo{Proof}
Assume that $PZ'>0$.
Since $(P+N)Z'=(f^*D+\Delta_x-Z)Z'=(\Delta_x-Z)Z'\leq 0$,
we have $NZ'<0$.
Thus $N$ contains $Z'$ as a component.
Since $PN=0$, this is contradiction.\qed
\enddemo\par
\proclaim{Lemma 3}
There exists an irreducible component $C$ of $f^{-1}(x)$ such that $PC>0$.
\endproclaim\par
\demo{Proof}
If $PC=0$ for all curves in $f^{-1}(x)$ then $P(Z-\Delta_x)=0$.
Hence $f^*D=P+(N+Z-\Delta_x)$ is also Zariski decomposition
of nef and big divisor $f^*D$.
Since $Z-\Delta_x$ is effective,
this contradicts the uniqueness of the Zariski decomposition.\qed
\enddemo\par
\remark{Remark}
This $C$ satisfies $(Z-\Delta_x)C<0$ by Lemma 2.
Moreover $(Z-\Delta_x)C=-1$ unless $x$ is of type $A_1$.
\endremark\par
First we consider the case of type $A_1$.
In this case, we have $\Delta_x=\bigl(1-\frac2{w_1}\bigr)C_1$
where $w_1=-C_1^2$.
Hence we have $(\Delta_x-Z)C_1=2$ and $\delta_x=4/w_1$.
Suppose $E_1C_1\geq 2$.
By Lemma 3, $N$ does not contain $C_1$ as a component.
Therefore we have $NC_1=(N+\Delta')C_1\geq E_1C_1\geq 2$.
On the other hand, we have $NC_1<(P+N)C_1=(\Delta_x-Z)C_1=2$,
which is absurd.
Since $E_1$ is integral near $Z$, we have $E_1C_1=1$.
Hence $f^*DE_1<(Z-\Delta_x)E_1=\frac2{w_1}=\frac12\delta_x$.
\par
Now, we assume $x$ is not of type $A_1$.
Thus the $C$ in Lemma 3 satisfies
$(\Delta_x-Z)C=1$ by the classification in the Appendix.
Hence there exists a component $C$ in $Z$
which satisfies $PC>0$ and $(\Delta_x-Z)C=1$.
\par
\proclaim{Lemma 4}
If $(Z-\Delta_x)E_1\geq\delta_x$ then there exists at least two components
in $Z$ which satisfies $(\Delta_x-Z)C_i=1$.
\endproclaim\par
\demo{Proof}
There exists a component $C$ which satisfies $PC>0$ and $(\Delta_x-Z)C=1$
by Lemma 3.
Hence $(Z-\Delta_x)P>0$ and
$(Z-\Delta_x)(N+\Delta')<(Z-\Delta_x)(P+N)=\delta_x$.
Since $E_1$ is taken from $N+\Delta'$,
there must be a component $C'$ in $N$ which satisfies $(Z-\Delta_x)C'<0$.
Hence this $C'$ is a component in $Z$ and $PC'=0$.
Thus these $C$ and $C'$ satisfies
$(\Delta_x-Z)C=(\Delta_x-Z)C'=1$.
Since $PC>0$ and $PC'=0$, $C$ and $C'$ are different.
\qed\enddemo\par
Let $T=\{i\mid PC_i>0,\text{ $C_i$ is a component of $Z$ }\}\ne\emptyset$.
We have $(\Delta_x-Z)C=1$ for any $i\in T$ by Lemma 2.
$T$ is a subset of $\{i\mid (\Delta_x-Z)C=1\}$,
which is indicated by $\bullet$ in the graph in the Appendix.
\par
\proclaim{Lemma 5}
Any component in $E_1$ does not meet $C_i$ for $i\in T$.
Thus we have $C_iE_1=0$ for $i\in T$.
\endproclaim
\demo{Proof}
Since $(P+N)C_i=(f^*D+\Delta_x-Z)C_i=1$ and $PC_i>0$,
we have $NC_i<1$.
Since $PC_i>0$, $N$ does not contain $C_i$ as a component.
Let $N=\sum n_jN_j+N'$ where $N_j$ is the components which meets $C_i$.
Then we have $n_j<1$ for all $j$.
Since $E_1$ is the divisor taken from $N+\Delta'$,
$E_1$ must consist of the components of $N'+\Delta'$.
Thus every component in $E_1$ does not meet $C_i$. \qed
\enddemo\par
Let $w_i=-C_i^2$.
Let $\mu_i=c_i-a_i$ where $c_i$ is the coefficient of $C_i$ in $Z$
and $a_i$ is the coefficient of $C_i$ in $\Delta_x$.
Since $PE_1=0$, we have
$$\align
  DE_2
    &= f^*DE_1 \leq (E_1+Z-\Delta_x)E_1 \\
    &= (E_1+\sum_{i\not\in T}\mu_i C_i+\sum_{i\in T}\mu_i C_i)E_1 \\
    &= \left( E_1+\frac{\sum_{i\not\in T}\mu_i C_i}{2}\right)^2
       -\frac14\left(\sum_{i\not\in T}\mu_i C_i\right)^2 \\
    &< -\frac14\left(\sum_{i\not\in T}\mu_i C_i\right)^2.
\endalign$$
Let $d=-\left(\sum_{i\not\in T}\mu_i C_i\right)^2$.
Hence it is enough to show $d\leq 2\delta_x$. \par
Case 1: $c_i=1$ for all $i\in T$.\par
Suppose $T=\{i,j,k\}$ where $i\ne j,\ j\ne k,\ k\ne i$.
By the classification in the Appendix,
this can occur only in case of Type $D_n$ or $E_n$,
and we have $C_iC_j=C_jC_k=C_kC_i=0$.
Let
$$\align
  d(t_1,t_2,t_3)
    &= -(Z-\Delta_x-t_1 C_i-t_2 C_j-t_3 C_k)^2 \\
    &= -(Z-\Delta_x)^2-2t_1+w_it_1^2-2t_2+w_jt_2^2-2t_3+w_kt_3^2.\\
\endalign$$
Since
$$
  -2=(K_X+C_i)C_i=(-\Delta_x+C_i)C_i < -(1-a_i)w_i,
$$
we have $0<\mu_i=1-a_i< \frac{2}{w_i}$.
Since $w_it_1^2-2t_1< 0$ when $0 < t_1 < \frac{2}{w_i}$,
we have $w_i\mu_i^2-2\mu_i < 0$.
The same inequality holds if we exchange $i$ to $j,k$.
Hence we have $d=d(\mu_i,\mu_j,\mu_k) < -(Z-\Delta)^2=\delta_x$.
\par
The cases of $\#T\leq 2$ are treated by the same method,
and we get $d\leq\delta_x$.\par
Note that if $x$ is of type $A_2$ and $T=\{1,2\}$ then
$\sum_{i\not\in T}\mu_i C_i=0$.\par
Case 2: $T=\{i,j\},\ i\ne j$ and $c_i=2$.\par
In this case, we have $c_j=1$ by the classification.
Let
$$\align
  d(t_1,t_2)
    &= -(Z-\Delta_x-t_1 C_i-t_2 C_j) \\
    &= -(Z-\Delta_x-t_1 C_i)^2-2t_2+w_jt_2^2. \\
\endalign$$
As in the case 1, we have $d=d(\mu_i,\mu_j) < -(Z-\Delta_x-\mu_i C_i)^2$.
So this case is treated as in the following case 3.
\par
Case 3: $T=\{i\}$ and $c_i=2$.\par
Case 3-1: There exists two components which satisfies $(\Delta_x-Z)C_i=1$
and $x$ is of type $D_n$.\par
Namely $x$ is of type of ($D_n$-ii) in the Appendix.\par
\proclaim{Lemma 6}
$a_1< \cdots < a_{j-1}\leq 2a_j-1$.
\endproclaim\par
\demo{Proof}
Since $\Delta_x C_j=-w_j+2$,
we have the equation $a_{j-1}-w_ja_j+a_{j+1}=-w_j+2$.
Since $a_j=a_{j+1}$ (cf\. proof of Proposition 1) and $w_j\geq 3$,
we have
$$
  a_{j-1}-a_j=2-w_j-2a_j+w_ja_j=(1-a_j)(2-w_j)\leq a_j-1.
$$
Hence we have $a_{j-1}\leq 2a_j-1<a_j$.\par
Similarly, since $\Delta_x C_i=-w_i+2$ for all $i\leq j-1$,
we have
$$\align
  a_{i-1}-a_i
    &= 2-w_i+w_ia_i-a_i-a_{i+1} \\
    &= a_i-a_{i+1}+(1-a_i)(2-w_i) \leq a_i-a_{i+1}.
\endalign$$
Since $a_{j-1}<a_j$, we have $a_{i-1}< a_i$ for all $i\leq j-1$.
Thus we have $a_1< \cdots < a_{j-1}\leq 2a_j-1$.\qed
\enddemo\par
Therefore we have
$$\align
  d &= -(Z-\Delta_x-(2-a_{j+1})C_{j+1})^2 \\
    &= -(Z-\Delta_x)^2-2(2-a_j)+2(2-a_j)^2
      \text{ since $a_j=a_{j+1}$ and $w_{j+1}=2$} \\
    &= 2-a_1+2(2-a_j)(1-a_j)
      \text{ since $\delta_x=2-a_1$} \\
    &= 2-a_1+\frac12(1-b)(3-b)
      \text{ where $b=2a_j-1$}.
\endalign$$
Let $f(t)=2-a_1+\frac12(1-t)(3-t)$.
Since $f$ is decreasing function on $[0,1]$,
we have $f(b)\leq f(a_1)$ by Lemma 6.
Hence
$$\align
  d &= f(b)\leq f(a_1) = 2-a_1+\frac12(1-a_1)(3-a_1) \\
    &= (2-a_1)\left(2-\frac12a_1\right)-\frac12 \\
    &< 2\delta_x.
\endalign$$
\par
Case 3-2: There exists only one component $C_i$
which satisfies $(\Delta_x-Z)C_i=1$ and $x$ is of type $D_n$.\par
Namely $x$ is of type ($D_n$-iii) in the Appendix.
In this case, we have $\delta_x=2-a_2$ (cf\. Prop\. 1) and
$$\align
  d &= -(Z-\Delta_x-(2-a_2)C_2)^2 \\
    &= -(Z-\Delta_x)^2-2(2-a_2)+2(2-a_2)^2 \\
    &= (2-a_2)(3-2a_2).
\endalign$$
If $x$ is not a rational double point, then $w_1\geq 3$.
Since $-w_1a_1+a_2=2-w_1$ and $a_1=a_2$, we have
$$
  a_1=a_2=1-\frac1{w_1-1}\geq\frac12.
$$
Therefore $d\leq 2(2-a_2)=2\delta_x$.\par
Suppose that $x$ is a rational double point.
Since the assumption of the case 3-2 and Lemma 4,
we have $ZE_1<2$.
Since $E_1$ is integral near $Z$, we have $ZE_1=1$.
Therefore
$$\align
  DE_2 &= f^*DE_1\leq (E_1+Z-\Delta_x)E_1 \\
       &= E_1^2+ZE_1 < 1 = \frac12\delta_x
\endalign$$
\par
Case 3-3: The dual graph is of type $E_n$.\par
We can calculate $d=-(Z-\Delta_x-\mu_iC_i)^2$ for $i\in T$ directly.
Then we have
$$\vbox{\offinterlineskip
  \halign{&\vrule#&\strut\;\hfil#\hfil\vbox{\hrule width0pt height 12pt}\cr
    \noalign{\hrule}
    & Type && 1 && 2 && 3 && 4 && 5 && 6 && 7 && 8
           && 9 && 10 && 11 && 12 && 13 && 14 && 15 & \cr
    height2pt & \omit && \omit && \omit && \omit && \omit
      && \omit && \omit && \omit && \omit && \omit && \omit
      && \omit && \omit && \omit && \omit && \omit & \cr
    \noalign{\hrule}
    & $\delta_x$ && $2$ && $\frac{14}{9}$ && $\frac65$ && $2$
      && $\frac85$ && $\frac{10}{7}$ && $\frac{12}{11}$ && $2$
      && $\frac{18}{11}$ && $\frac{23}{13}$ && $\frac{33}{23}$
      && $\frac{11}{7}$ && $\frac{21}{17}$ && $\frac{26}{19}$
      && $\frac{30}{29}$ & \cr
    height2pt & \omit && \omit && \omit && \omit && \omit
      && \omit && \omit && \omit && \omit && \omit && \omit
      && \omit && \omit && \omit && \omit && \omit & \cr
    \noalign{\hrule}
    & $\frac{d}4$ && $1.5$ && $\frac{7.3\ccdot}{9}$ && $\frac{2.1}{5}$
      && $1.5$ && $\frac{4.4}{5}$ && $\frac{4.6\ccdot}{7}$
      && $\frac{3.5\ccdot}{11}$ && $1.5$ && $\frac{10.2\ccdot}{11}$
      && $\frac{9.2\ccdot}{13}$ && $\frac{8.7\ccdot}{23}$
      && $\frac{5.8\ccdot}{7}$ && $\frac{7.7\ccdot}{17}$
      && $\frac{11.2\ccdot}{19}$ && $\frac{8.0\ccdot}{29}$ & \cr
    height2pt & \omit && \omit && \omit && \omit && \omit
      && \omit && \omit && \omit && \omit && \omit && \omit
      && \omit && \omit && \omit && \omit && \omit & \cr
    \noalign{\hrule}
}}$$
If $x$ is a rational double point then
we have $DE_2<1$ as in the case 3-2.
Hence if $x$ is not neither of the type 2, 5, 9 or 12 then $DE_2<\delta_x/2$.
\par
Thus we consider the case of the type 2, 5, 9 and 12.
Since $DE_2<(Z-\Delta_x)E_1$,
if $(Z-\Delta_x)E_1\leq\delta_x/2$ then there is nothing to show.
Hence we consider the case $\delta_x/2<(Z-\Delta_x)E_1<\delta_x$
by Lemma 4.
Let $e=(Z-\Delta_x)E_1$.\par
We say that $Z-\Delta_x$ is
$(\mu_0;\mu_1,\cdots,\mu_{n_1};\mu_1',\cdots,\mu_{n_2}';\mu_1'')$
if
$$
  Z-\Delta_x=\mu_0 C_0+\sum_{i=1}^{n_1} \mu_iC_i+
    \sum_{j=1}^{n_2} \mu_j'C_j'+\mu_1''C_1''.$$
\par
The case of type 2:\par
In this case, $\delta_x=\frac{14}{9}$, $T=\{2\}$ and $Z-\Delta_x$ is
$$
  \left(\frac43;\frac79,\frac{14}{9};\frac49;\frac23\right).
$$
Let $E_1=\sum b_jC_j+\sum e_iE_i'+E''$ where $E_i'$ are the components
which meet $Z$ but not exceptional, and $E''$ is disjoint from $Z$.
Set also $\sum e_iE_i'+E''=E'$.
Since $E_1$ does not contain the component which satisfies
$(Z-\Delta_x)C_i=-1$,
we have $e=(Z-\Delta_x)E_1=(Z-\Delta_x)E'$.
Since $E_1$ is integral near $Z$, $e$ must be the sum of some coefficients of
$Z-\Delta_x$.
Hence we consider only the case $9e=8,9,10,11,12,13$.
If $9e=13$ and $E_1$ has a component which meets $C_1$ and
a component which meets $C_1''$ then
we write ``13=7+6'' by using the coefficients of such components.\par
Suppose that $E_1\cap C_1=\emptyset$.
By Lemma 5, we have $C_1E_1=C_2E_1=0$.
Thus
$$\align
  f^*DE_1 &\leq (E_1+Z-\Delta_x)E_1 \\
    &= \left(E_1+\frac43C_0+\frac49C_1'+\frac23C_1''\right)E_1 \\
    &< -\frac14\left(\frac43C_0+\frac49C_1'+\frac23C_1''\right)^2 \\
    &= \frac{4.66\cdots}{9}.
\endalign$$
Hence we may assume that there is a component which meet $C_1$.
This occurs only in the case of ``13=7+6''
namely $E'C_1=E'C_1''=1,\ E'C_{\bullet}^{\bullet}=0$ for other
$C_{\bullet}^{\bullet}$ where ${}_{\bullet}^{\bullet}$ means index and primes.
or ``11=7+4.'' namely $E'C_1=E'C_1'=1,\ E'C_{\bullet}^{\bullet}=0$ for others.
Let $E_1=b_1'C_1'+b_1''C_1''+E'$ where $E'$ consists of the components
which is not contained in $Z$, since $E_1$ does not contain
the components which meets $C_2$ by Lemma 5.
Thus we have $C_1'E'\leq 1$ and $C_1''E'\leq 1$.
Hence we have
$$\align
  f^*DE_1 &\leq (E_1+Z-\Delta_x)E_1 \\
    &= (E'+b_1'C_1'+b_1''C_1''+Z-\Delta_x)(E'+b_1'C_1'+b_1''C_1'') \\
    &= (E'+Z-\Delta_x)E'+2E'(b_1'C_1'+b_1''C_1'')+(b_1'C_1'+b_1''C_1'')^2 \\
    &\leq (E'+Z-\Delta_x)E'+2b_1'-3(b_1')^2+2b_1''-2(b_1'')^2 \\
    &\leq (E'+Z-\Delta_x)E'\text{ since $b_1'$ and $b_1''$ are integral } \\
    &= \left(E'+\frac79C_1+\frac49C_1'+\frac23C_1''\right)E' \\
    &< -\frac14\left(\frac79C_1+\frac49C_1'+\frac23C_1''\right)^2 \\
    &= \frac{6.05\cdots}{9}.
\endalign$$
Therefore $f^*DE_1<\delta_x/2$.\par
The case of type 5:\par
In this case, $\delta_x=\frac85$, $T=\{2\}$ and $Z-\Delta_x$ is
$$
  \left(\frac65;\frac45,\frac85,\frac75;\frac25;\frac35\right).
$$
Suppose that $E_1\cap C_1=\emptyset$.
By Lemma 5, we have $C_1E_1=C_2E_1=0$.
Thus
$$\align
  f^*DE_1 &\leq (E_1+Z-\Delta_x)E_1 \\
    &= \left(E_1+\frac75C_3+\frac65C_0+\frac25C_1'+\frac35C_1''\right)E_1 \\
    &< -\frac14\left(\frac75C_3+\frac65C_0+\frac25C_1'+\frac35C_1''\right)^2 \\
    &= \frac{2.8}{5}.
\endalign$$
Hence we may assume that $C_1E_1>0$.
This occurs only in the case of ``7=4+3'' or ``6=4+2.''
Let $E_1=b_0C_0+b_1'C_1'+b_1''C_1''+E'$ where $E'$ is as above.
Thus we have $C_1'E'\leq 1$ and $C_1''E'\leq 1$.\par
If $b_0=0$ then
$$\align
  f^*DE_1 &\leq (E'+b_1'C_1'+b_1''C_1''+Z-\Delta_x)(E'+b_1'C_1'+b_1''C_1'') \\
    &\leq (E'+Z-\Delta_x)E'+2b_1'-3(b_1')^2+2b_1''-2(b_1'')^2 \\
    &\leq \left(E'+\frac45C_1+\frac25C_1'+\frac35C_1''\right)E' \\
    &< -\frac14\left(\frac45C_1+\frac25C_1'+\frac35C_1''\right)^2 \\
    &= \frac{3.1}{5}.
\endalign$$\par
If $b_0\ne 0$ then
$$\align
  \left(\frac65C_0+\frac25C_1'+\frac35C_1''\right)E_1
    &= \left(\frac65C_0+\frac25C_1'+\frac35C_1''\right)
         (E'+b_0C_0+b_1'C_1'+b_1''C_1'') \\
    &= E'(\frac25C_1'+\frac35C_1'')-\frac75b_0
    \leq -\frac45.
\endalign$$
Therefore
$$\align
  f^*DE_1
    &\leq \left(E_1+\frac45C_1+\frac75C_3+\frac65C_0
      +\frac25C_1'+\frac35C_1''\right)E_1 \\
    &\leq \left(E_1+\frac45C_1+\frac75C_3\right)E_1-\frac45 \\
    &< -\frac14\left(\frac45C_1+\frac75C_3\right)^2-\frac45 \\
    &= \frac{2.5}{5}.
\endalign$$
Therefore $f^*DE_1<\delta_x/2$.\par
The case of type 9:\par\nopagebreak
In this case, $\delta_x=\frac{18}{11}$, $T=\{2\}$ and $Z-\Delta_x$ is
$$
  \left(\frac{12}{11};\frac9{11},\frac{18}{11},\frac{16}{11},\frac{14}{11};
    \frac4{11};\frac6{11}\right).
$$
Suppose that $E_1\cap C_1=\emptyset$.
By Lemma 5, we have $C_1E_1=C_2E_1=0$.
Thus
$$\align
  f^*DE_1 &\leq (E_1+Z-\Delta_x)E_1 \\
    &= \left(E_1+\frac{16}{11}C_3+\frac{14}{11}C_4+\frac{12}{11}C_0
         +\frac{4}{11}C_1'+\frac{6}{11}C_1''\right)E_1 \\
    &< -\frac14\left(\frac{16}{11}C_3+\frac{14}{11}C_4+\frac{12}{11}C_0
         +\frac{4}{11}C_1'+\frac{6}{11}C_1''\right)^2 \\
    &= \frac{6.54\cdots}{11}.
\endalign$$
We assume that there is a component which meet $C_1$.
This occurs only in the case of ``15=9+6'', ``17=9+4+4'' or ``13=9+4.''
Let $E_1=b_4C_4+b_0C_0+b_1'C_1'+b_1''C_1''+E'$ where $E'$ is as above.
Thus we have $C_1'E'\leq 2$ and $C_1''E'\leq 1$.\par
First we assume $b_0=0$.
In the case of ``15=9+6'' or ``13=9+4'' we have
$$\align
  f^*DE_1
    &\leq (E'+b_4C_4+b_1'C_1'+b_1''C_1''+Z-\Delta_x)
      (E'+b_4C_4+b_1'C_1'+b_1''C_1'') \\
    &\leq (E'+Z-\Delta_x)E'+2b_1'-3(b_1')^2+2b_1''-2(b_1'')^2-2(b_4)^2 \\
    &\leq (E'+Z-\Delta_x)E' \\
    &= \left(E'+\frac9{11}C_1+\frac4{11}C_1'+\frac6{11}C_1''\right)E' \\
    &< -\frac14\left(\frac9{11}C_1+\frac4{11}C_1'+\frac6{11}C_1''\right)^2 \\
    &= \frac{6.40\cdots}{11}
\endalign$$
Suppose $E_1$ is of the case of ``17=9+4+4.''
Since $P(E'+t_1C_1+t_3C_3+t_4C_4+t_0C_0+t_1'C_1'+t_1''C_1'')=0$
for any $t_1,t_3,t_4,t_0,t_1',t_1''$,
we have
$$\align
  0 &> (E' + t_1C_1+t_3C_3+t_4C_4+t_0C_0+t_1'C_1'+t_1''C_1'')^2 \\
    &={E'}^2+2t_1+4t_1'-2t_1^2-2t_3^2-2t_4^2-2t_0^2-3{t_1'}^2-2{t_1''}^2 \\
    &\quad +2t_3t_4+2t_4t_0+2t_0t_1'+2t_0t_1''.
\endalign$$
If we take $(t_1,t_3,t_4,t_0,t_1',t_1'')=(1/2,4/9,8/9,4/3,10/9,2/3)$,
we get ${E'}^2+\frac{49}{18}<0$.
Hence
$$\align
  f^*DE_1
    &\leq (E'+b_4C_4+b_1'C_1'+b_1''C_1''+Z-\Delta_x)
      (E'+b_4C_4+b_1'C_1'+b_1''C_1'') \\
    &\leq (E'+Z-\Delta_x)E'+4b_1'-3(b_1')^2+2b_1''-2(b_1'')^2-2(b_4)^2 \\
    &\leq (E'+Z-\Delta_x)E'+1 \\
    &= (E'+\frac9{11}C_1+\frac4{11}C_1')E'+1 \\
    &= {E'}^2+\frac{17}{11}+1 \\
    &< -\frac{49}{18}+\frac{17}{11}+1 = -\frac{35}{198}
\endalign$$
which is absurd.\par
We assume that $b_0\ne 0$.
If $b_4=0$ then
$$
  \left(\frac{12}{11}C_0+\frac4{11}C_1'+\frac6{11}C_1''\right)E_1
    = E'\left(\frac4{11}C_1'+\frac6{11}C_1''\right)-\frac{14}{11}b_0
    \leq -\frac6{11}
$$
Thus we have
$$\align
  f^*DE_1
    &\leq \left(E_1+\frac9{11}C_1+\frac{16}{11}C_3+\frac{14}{11}C_4
       +\frac{12}{11}C_0+\frac4{11}C_1'+\frac6{11}C_1''\right)E_1 \\
    &\leq \left(E_1+\frac9{11}C_1+\frac{16}{11}C_3+\frac{14}{11}C_4\right)E_1
       -\frac6{11} \\
    &< -\frac14\left(\frac9{11}C_1+\frac{16}{11}C_3+\frac{14}{11}C_4\right)^2
       -\frac6{11} \\
    &= \frac{8.04\cdots}{11}.
\endalign$$
If $b_4\ne 0$ then
$$
  \left(\frac{14}{11}C_4+\frac{12}{11}C_0+\frac4{11}C_1'
    +\frac6{11}C_1''\right)E_1
  = E'\left(\frac4{11}C_1'+\frac6{11}C_1''\right)-\frac{16}{11}b_4
  \leq -\frac8{11}.
$$
Thus we have
$$\align
  f^*DE_1
    &\leq \left(E_1+\frac9{11}C_1+\frac{16}{11}C_3+\frac{14}{11}C_4
      +\frac{12}{11}C_0+\frac4{11}C_1'+\frac6{11}C_1''\right)E_1 \\
    &\leq \left(E_1+\frac9{11}C_1+\frac{16}{11}C_3\right)E_1-\frac8{11} \\
    &< -\frac14\left(\frac9{11}C_1+\frac{16}{11}C_3\right)^2-\frac8{11} \\
    &= \frac{7.31\cdots}{11}.
\endalign$$
Therefore $f^*DE_1<\delta_x/2$.\par
The case of type 12:\par
In this case, $\delta_x=\frac{11}7$, $T=\{1''\}$ and $Z-\Delta_x$ is
$$
  \left(\frac{15}7;\frac37,\frac97;\frac57,\frac{10}7;\frac{11}7\right).
$$
Let $E_1=b_1C_1+b_2C_2+b_1'C_1'+b_2'C_2'+E'$ where $E'$ is as above. \par
First, we consider the case of ``8=5+3.''\par
If $b_2\ne 0$ then
$$\align
   \left(\frac37C_1+\frac97C_2\right)E_1
     &= \left(\frac37C_1+\frac97C_2\right)E'-\frac{15}7b_2 \\
     &= \frac37-\frac{15}7b_2
     \leq -\frac{12}7.
\endalign$$
Hence
$$\align
  f^*DE_1 &\leq (E_1+Z-\Delta_x)E_1 \\
    &\leq \left(E_1+\frac{15}7C_0+\frac{10}7C_2'+\frac57C_1'
      \right)E_1-\frac{12}7 \\
    &< -\frac14\left(\frac{15}7C_0+\frac{10}7C_2'+\frac57C_1'
      \right)^2-\frac{12}7 \\
    &= -\frac9{49} <0.
\endalign$$
That is contradiction.
Thus $b_2=0$.\par
If $b_2'\ne 0$ then
$$\align
  \left(\frac{10}7C_2'+\frac57C_1'\right)E_1
    &= \left(\frac{10}7C_2'+\frac57C_1'\right)E'-\frac{15}7b_2' \\
    &= \frac57-\frac{15}7b_2'
    \leq -\frac{10}7.
\endalign$$
Hence
$$\align
  f^*DE_1 &\leq (E_1+Z-\Delta_x)E_1 \\
    &\leq \left(E_1+\frac37C_1+\frac97C_2+\frac{15}7C_0
      \right)E_1-\frac{10}7 \\
    &< -\frac14\left(\frac37C_1+\frac97C_2+\frac{15}7C_0
      \right)^2-\frac{10}7 \\
    &= \frac{1.25}7.
\endalign$$
If $b_2=b_2'=0$ then
$$\align
  f^*DE_1 &\leq (E_1+Z-\Delta_x)E_1 \\
    &= (E'+b_1C_1+b_1'C_1'+Z-\Delta_x)(E'+b_1C_1+b_1'C_1') \\
    &= (E'+Z-\Delta_x)E'+2E'(b_1C_1+b_1'C_1')+(b_1C_1+b_1'C_1')^2 \\
    &= (E'+Z-\Delta_x)E'+2b_1-3b_1^2+2b_1'-2b_1^2 \\
    &\leq (E'+Z-\Delta_x)E' \\
    &= \left(E'+\frac37C_1+\frac57C_1'\right)E' \\
    &< -\frac14\left(\frac37C_1+\frac57C_1'\right)^2 \\
    &= \frac{2.75}7.
\endalign$$
\par
Next, we consider the cases of ``9=9'', ``9=3+3+3'' or ``6=3+3.''\par
If $b_2'\ne 0$ then
$$
   \left(\frac{10}7C_2'+\frac57C_1'\right)E_1
     = \left(\frac{10}7C_2'+\frac57C_1'\right)E'-\frac{15}7b_2'
     \leq -\frac{15}7.
$$
Hence
$$\align
  f^*DE_1 &\leq (E_1+Z-\Delta_x)E_1 \\
    &\leq \left(E_1+\frac37C_1+\frac97C_2+\frac{15}7C_0\right)E_1-\frac{15}7 \\
    &< -\frac14\left(\frac37C_1+\frac97C_2+\frac{15}7C_0\right)^2-\frac{15}7 \\
    &= \frac{45}{28}-\frac{15}7 <0.
\endalign$$
That is contradiction.
Thus $b_2'=0$.\par
If $b_2\ne 0$ then
$$\align
  \left(\frac37C_1+\frac97C_2\right)E_1
    &= \left(\frac37C_1+\frac97C_2\right)E'-\frac{15}7b_2 \\
    &\leq \frac97-\frac{15}7b_2'
    \leq -\frac67.
\endalign$$
Hence
$$\align
  f^*DE_1 &\leq (E_1+Z-\Delta_x)E_1 \\
    &\leq \left(E_1+\frac{15}7C_0+\frac{10}7C_2'+\frac57C_1'\right)E_1-\frac67
\\
    &< -\frac14\left(\frac{15}7C_0+\frac{10}7C_2'+\frac57C_1'\right)^2-\frac67
\\
    &= \frac{4.71\cdots}7.
\endalign$$
If $b_2=b_2'=0$ then
$$\align
  f^*DE_1 &\leq (E_1+Z-\Delta_x)E_1 \\
    &= (E'+b_1C_1+b_1'C_1'+Z-\Delta_x)(E'+b_1C_1+b_1'C_1') \\
    &= (E'+b_1C_1+Z-\Delta_x)(E'+b_1C_1)+2(E'+b_1C_1)(b_1'C_1')+(b_1'C_1')^2 \\
    &\leq (E'+b_1C_1+Z-\Delta_x)(E'+b_1C_1) \\
    &= \left((E'+b_1C_1)+\left(\frac37C_1+\frac97C_2\right)\right)(E'+b_1C_1)
\\
    &< -\frac14\left(\frac37C_1+\frac97C_2\right)^2 \\
    &= \frac{4.82\cdots}7.
\endalign$$
\par
Finally, we consider the cases of ``10=10'' or ``10=5+5.''\par
If $b_2\ne 0$ then
$$
   \left(\frac37C_1+\frac97C_2\right)E_1
     = \left(\frac37C_1+\frac97C_2\right)E'-\frac{15}7b_2
     \leq -\frac{15}7.
$$
Hence
$$\align
  f^*DE_1 &\leq (E_1+Z-\Delta_x)E_1 \\
    &\leq \left(E_1+\frac{15}7C_0+\frac{10}7C_2'+\frac57C_1'\right)E_1
      -\frac{15}7 \\
    &< -\frac14\left(\frac{15}7C_0+\frac{10}7C_2'+\frac57C_1'\right)^2
      -\frac{15}7 \\
    &= \frac{75}{49}-\frac{15}7 <0.
\endalign$$
That is contradiction.
Thus $b_2=0$.\par
If $b_2'=0$ then
$$\align
  f^*DE_1 &\leq (E_1+Z-\Delta_x)E_1 \\
    &= (E'+b_1C_1+b_1'C_1'+Z-\Delta_x)(E'+b_1C_1+b_1'C_1') \\
    &= (E'+b_1'C_1'+Z-\Delta_x)(E'+b_1'C_1')+2(E'+b_1'C_1')(b_1C_1)+(b_1C_1)^2
\\
    &\leq (E'+b_1'C_1'+Z-\Delta_x)(E'+b_1'C_1') \\
    &= \left((E'+b_1'C_1')+\left(\frac{10}7C_2'+\frac57C_1'\right)\right)
      (E'+b_1'C_1') \\
    &< -\frac14\left(\frac{10}7C_2'+\frac57C_1'\right)^2 \\
    &= \frac{5.35\cdots}7.
\endalign$$
Hence we assume $b_2'\ne 0$.
Thus we have
$$
  \left(\frac{10}7C_2'+\frac57C_1'\right)E_1
     = \left(\frac{10}7C_2'+\frac57C_1'\right)E'-\frac{15}7b_2'
     \leq -\frac57.
$$
If $b_1\ne 0$ then
$$
  \frac37C_1E_1=-\frac97b_1\leq-\frac97.
$$
Therefore
$$\align
  f^*DE_1 &\leq (E_1+Z-\Delta_x)E_1 \\
    &\leq \left(E_1+\frac97C_2+\frac{15}7C_0\right)E_1-2 \\
    &< -\frac14\left(\frac97C_2+\frac{15}7C_0\right)^2-2 \\
    &= \frac{171}{98}-2 < 0.
\endalign$$
That is contradiction.
Thus $b_1=0$.
Furthermore if $b_2'\geq 2$ then
$$
  \left(\frac{10}7C_2'+\frac57C_1'\right)E_1 \leq -\frac{20}7.
$$
Hence we have
$$\align
  f^*DE_1 &\leq (E_1+Z-\Delta_x)E_1 \\
    &\leq \left(E_1+\frac37C_1+\frac97C_2+\frac{15}7C_0\right)E_1-\frac{20}7 \\
    &< -\frac14\left(\frac37C_1+\frac97C_2+\frac{15}7C_0\right)^2-\frac{20}7 \\
    &= \frac{11.25}7-\frac{20}7 < 0
\endalign$$
That is contradiction.
Thus we assume $b_1=b_2=0$, $b_2'=1$ and $E_1=E'+b_1'C_1'+C_2'$.
In case ``10=10'' we have
$$\align
  (E' &+ b_1'C_1'+C_2'+Z-\Delta_x)(E'+b_1'C_1'+C_2') \\
    &= (E'+Z-\Delta_x)E'+2E'(b_1'C_1'+C_2')+(b_1'C_1+C_2')^2 \\
    &= (E'+Z-\Delta_x)E'+2b_1'-2{b_1'}^2 \\
    &\leq \left(E'+\frac{10}7C_2'+\frac57C_1'\right)E' \\
    &< -\frac14\left(\frac{10}7C_2'+\frac57C_1'\right)^2 \\
    &= \frac{5.35\cdots}7.
\endalign$$
Hence we consider the case of ``10=5+5.''
Since
$$
  \frac57C_1'E_1=\frac57C_1'(E'+b_1'C_1'+C_2')
    =\frac{15}7-\frac{10}7b_1'
$$
if $b_1'\geq 2$ then
$$\align
  f^*DE_1 &\leq (E_1+Z-\Delta_x)E_1 \\
    &\leq
\left(E_1+\frac37C_1+\frac97C_2+\frac{15}7C_0+\frac{10}7C_2'\right)E_1
      -\frac57 \\
    &<
-\frac14\left(\frac37C_1+\frac97C_2+\frac{15}7C_0+\frac{10}7C_2'\right)^2
      -\frac57 \\
    &= \frac{2.67\cdots}7.
\endalign$$
Hence we assume $b_1'\leq 1$.
In this case, we have
$$\align
  (E' &+ b_1'C_1'+C_2'+Z-\Delta_x)(E'+b_1'C_1'+C_2') \\
    &= (E'+b_1'C_1'+Z-\Delta_x)(E'+b_1'C_1')+2(E'+b_1'C_1')C_2'+(C_2')^2 \\
    &= (E'+b_1'C_1'+Z-\Delta_x)(E'+b_1'C_1')+2b_1'-2 \\
    &\leq \left(E'+b_1'C_1'+\frac{10}7C_2'+\frac57C_1'\right)(E'+b_1'C_1') \\
    &< -\frac14\left(\frac{10}7C_2'+\frac57C_1'\right)^2 \\
    &= \frac{5.35\cdots}7.
\endalign$$
Hence we have $f^*DE_1<\delta_x/2$.\par
Therefore, for all cases we have $DE_2<\delta_x/2$.
If we replace $E_2$ by $E$, then we get the following bound.
$$\gather
 0\leq DE<\frac12\delta_x \\
 DE-\frac14\delta_x\leq E^2\leq\frac{(DE)^2}{D^2} \\
 E^2<0 \text{ if } DE=0.
\endgather$$
Thus we complete the proof when $x$ is log terminal singular point.
Hence we complete the proof of the theorem.\qed
\par\medskip
{}From this theorem, we can get the following results immediately.\par
\proclaim{Corollary 1}
Let $Y$ be a normal surface over $k$.
Let $A$ be an ample divisor on $Y$.
If $t\geq 3$ and $K_Y+tA$ is a Cartier divisor
then $|K_Y+tA|$ is free.
Moreover, if $A^2>1$ and $K_Y+2A$ is a Cartier divisor
then $|K_Y+2A|$ is free.
\endproclaim\par
\proclaim{Corollary 2}
Let $Y$ be a normal surface over $k$.
Let $D$ be a nef divisor on $Y$.
Fix one point $x$ in $Y$.
If $D$ satisfies
\roster
\item $D^2>\delta_x$
\item $DC\geq \frac12\delta_x$ for all irreducible curves $C$ on $Y$
\item $K_Y+D$ is a Cartier divisor
\endroster
then $|K_Y+D|$ is free at $x$.
\endproclaim\par
\proclaim{Corollary 3}
Let $Y$ be a normal surface having only rational double points as
singularities.
If $D$ is an ample Cartier divisor and $D^2>2$ then
every base points of $|K_Y+D|$ is smooth.
\endproclaim\par
\proclaim{Corollary 4}
Let $Y$ be a normal surface over $k$.
Let $D$ be a nef divisor on $Y$ such that $K_Y+D$ is a Cartier divisor.
Let $D^2>2$.
If $x$ is a rational double point of type $E_8$ then
$|K_Y+D|$ is free at $x$.
\endproclaim\par
\demo{Proof}
Suppose that $x$ is a base point of $|K_Y+D|$.
Since $E_1$ in theorem is integral near $Z$ and
all coefficients of the fundamental cycle of type $E_8$
are greater than or equal $2$.
Thus we have $ZE_1\geq 2$.
Hence by Lemma 4, there must be at least two components in $Z$
which satisfies $-ZC=1$, which is contradiction.
Thus $|K_Y+D|$ is free at $x$.\qed
\enddemo\par
\newpage
\head Appendix \endhead
\centerline{Table 1}
$$\vbox{\offinterlineskip
  \halign{\vrule#&\strut\ \hfil#\hfil &\vrule#&\ #\hfil&\vrule#&
          \ \hfil#\hfil &\vrule#&
          \ \hfil#\hfil\vbox{\hrule width0pt height 12pt}&\vrule# \cr
    \noalign{\hrule}
    & type && dual graph && $\delta_x$ && $\mu=3$ & \cr
    height2pt & \omit && \omit && \omit && \omit & \cr
    \noalign{\hrule}
    & 1 && $(\mu;2,2;2,2)$ && $\frac{11}{6}+\frac{1}{6(6\mu-11)}$
        && $\frac{13}{7}$ & \cr
    height2pt & \omit && \omit && \omit && \omit & \cr
    \noalign{\hrule}
    & 2 && $(\mu;2,2;3)$ && $\frac{3}{2}-\frac{1}{18(2\mu-3)}$
        && $\frac{41}{27}$ & \cr
    height2pt & \omit && \omit && \omit && \omit & \cr
    \noalign{\hrule}
    & 3 && $(\mu;3;3)$ && $\frac{7}{6}+\frac{1}{6(6\mu-7)}$
        && $\frac{13}{11}$ & \cr
    height2pt & \omit && \omit && \omit && \omit & \cr
    \noalign{\hrule}
    & 4 && $(\mu;2,2,2;2,2)$ && $\frac{23}{12}+\frac{1}{12(12\mu-23)}$
        && $\frac{25}{13}$ & \cr
    height2pt & \omit && \omit && \omit && \omit & \cr
    \noalign{\hrule}
    & 5 && $(\mu;2,2,2;3)$ && $\frac{19}{12}+\frac{1}{12(12\mu-19)}$
        && $\frac{27}{17}$ & \cr
    height2pt & \omit && \omit && \omit && \omit & \cr
    \noalign{\hrule}
    & 6 && $(\mu;4;2,2)$ && $\frac{17}{12}+\frac{1}{12(12\mu-17)}$
        && $\frac{27}{19}$ & \cr
    height2pt & \omit && \omit && \omit && \omit & \cr
    \noalign{\hrule}
    & 7 && $(\mu;4;3)$ && $\frac{13}{12}+\frac{1}{12(12\mu-13)}$
        && $\frac{25}{23}$ & \cr
    height2pt & \omit && \omit && \omit && \omit & \cr
    \noalign{\hrule}
    & 8 && $(\mu;2,2,2,2;2,2)$ && $\frac{59}{30}+\frac{1}{30(30\mu-59)}$
        && $\frac{61}{31}$ & \cr
    height2pt & \omit && \omit && \omit && \omit & \cr
    \noalign{\hrule}
    & 9 && $(\mu;2,2,2,2;3)$ && $\frac{49}{30}+\frac{1}{30(30\mu-49)}$
        && $\frac{67}{41}$ & \cr
    height2pt & \omit && \omit && \omit && \omit & \cr
    \noalign{\hrule}
    & 10 && $(\mu;2,3;2,2)$ && $\frac{53}{30}+\frac{1}{30(30\mu-47)}$
         && $\frac{76}{43}$ & \cr
    height2pt & \omit && \omit && \omit && \omit & \cr
    \noalign{\hrule}
    & 11 && $(\mu;2,3;3)$ && $\frac{43}{30}+\frac{1}{30(30\mu-37)}$
         && $\frac{76}{53}$ & \cr
    height2pt & \omit && \omit && \omit && \omit & \cr
    \noalign{\hrule}
    & 12 && $(\mu;3,2;2,2)$ && $\frac{47}{30}+\frac{1}{30(30\mu-53)}$
         && $\frac{58}{37}$ & \cr
    height2pt & \omit && \omit && \omit && \omit & \cr
    \noalign{\hrule}
    & 13 && $(\mu;3,2;3)$ && $\frac{37}{30}+\frac{1}{30(30\mu-43)}$
         && $\frac{58}{47}$ & \cr
    height2pt & \omit && \omit && \omit && \omit & \cr
    \noalign{\hrule}
    & 14 && $(\mu;5;2,2)$ && $\frac{41}{30}+\frac{1}{30(30\mu-41)}$
         && $\frac{67}{49}$ & \cr
    height2pt & \omit && \omit && \omit && \omit & \cr
    \noalign{\hrule}
    & 15 && $(\mu;5;3)$ && $\frac{31}{30}+\frac{1}{30(30\mu-31)}$
         && $\frac{61}{59}$ & \cr
    height2pt & \omit && \omit && \omit && \omit & \cr
    \noalign{\hrule}
    }}
$$
$$
  \delta_x=-(\Delta-Z)^2,\quad\text{if } \mu\geq 3.
$$
\par
\newpage
\topinsert
\vskip 9truein
\includegraphics{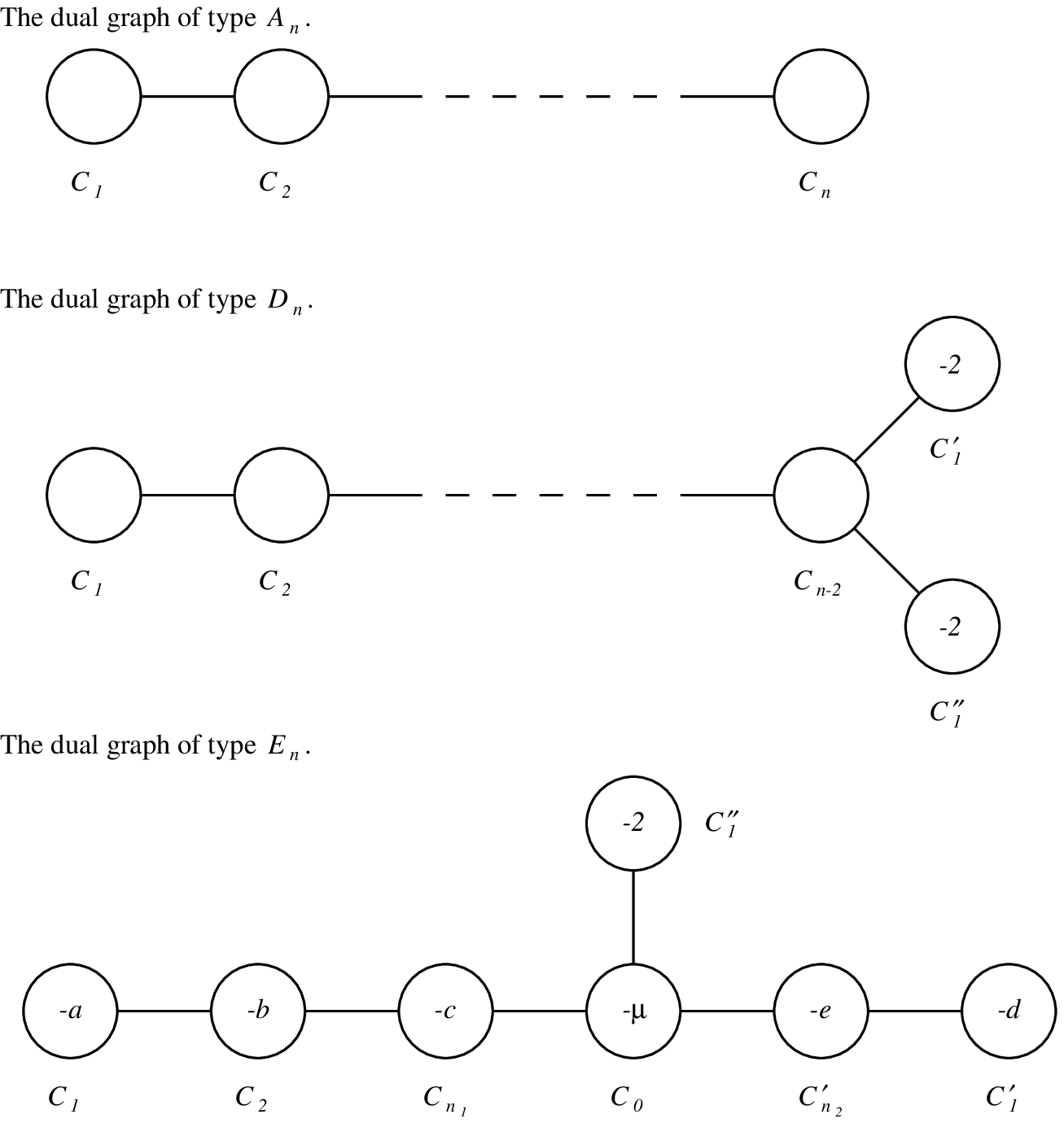}
\endinsert
\newpage
\topinsert
\vskip 9truein
\includegraphics{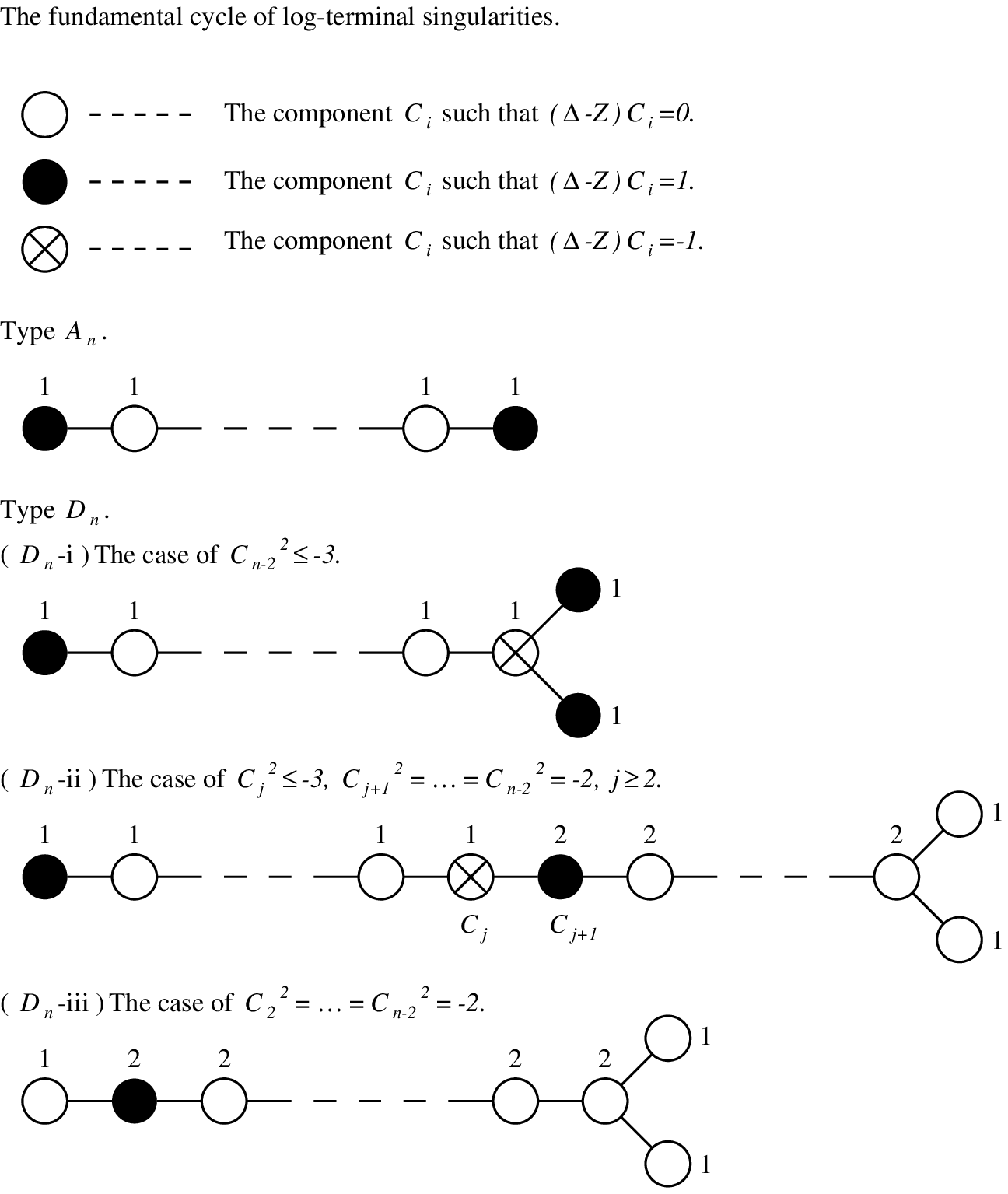}
\endinsert
\newpage
\topinsert
\vskip 9truein
\includegraphics{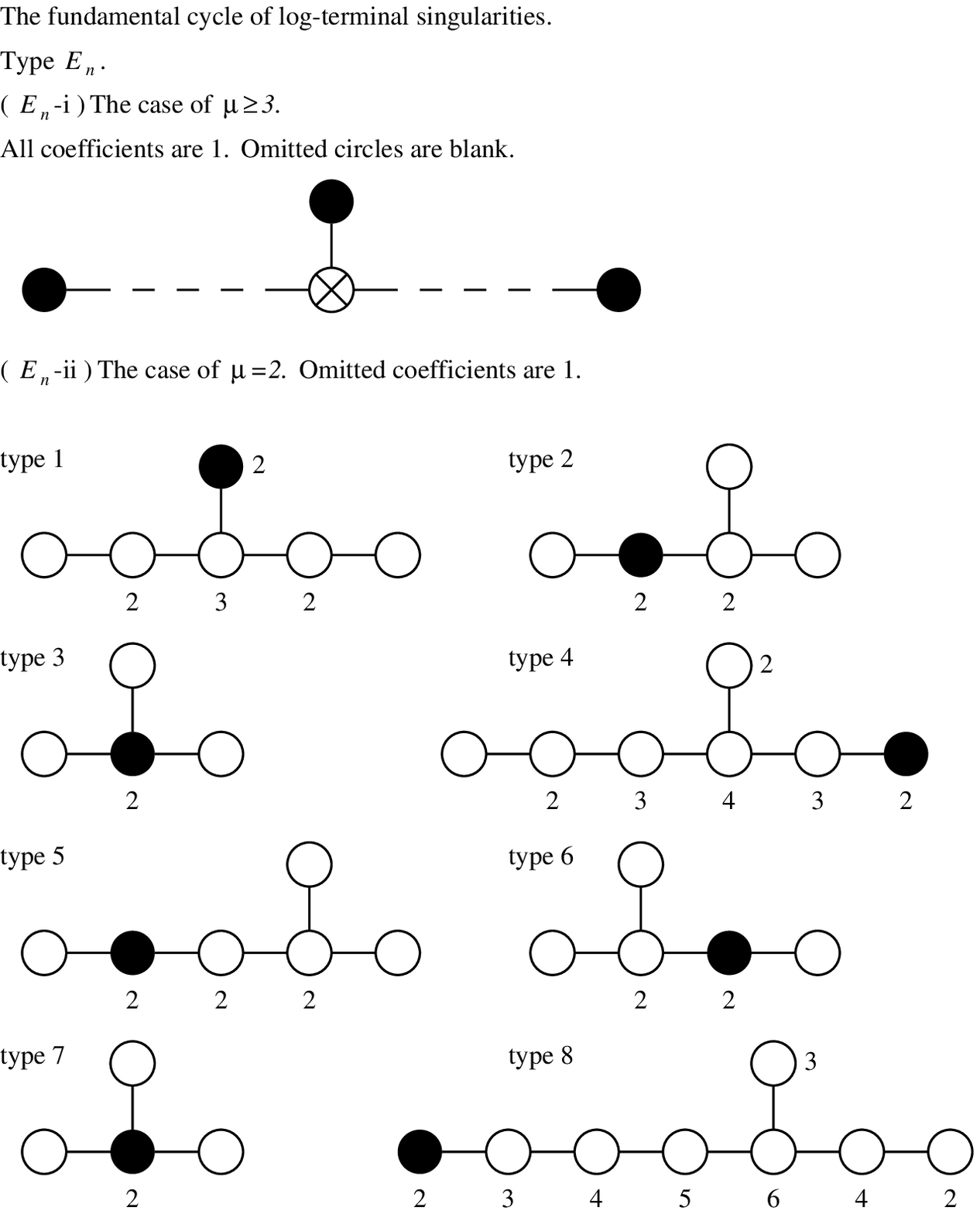}
\endinsert
\newpage
\topinsert
\vskip 9truein
\includegraphics{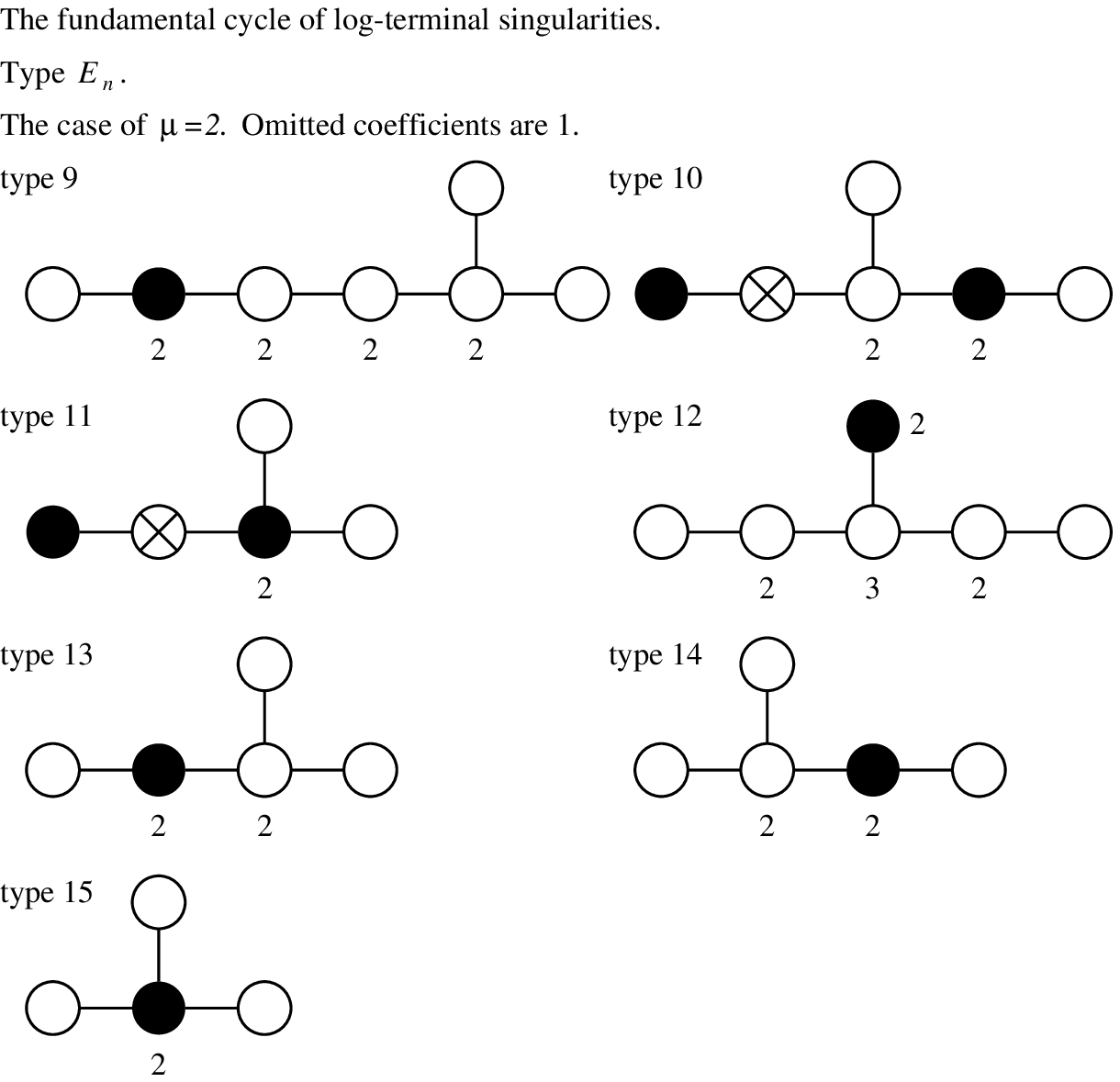}
\endinsert
\newpage

\enddocument